\documentclass[a4paper,DIV=14]{scrartcl}
\usepackage[T1]{fontenc}
\usepackage[utf8]{inputenc}
\usepackage[english]{babel}
\usepackage{url}
\usepackage[
  algo2e,
  lined,
  boxed,
  commentsnumbered,
  ruled
]{algorithm2e}
\usepackage{amsmath}
\usepackage{amssymb}
\usepackage{amsthm}
\theoremstyle{plain}
\newtheorem{theorem}{Theorem}[section]
\newtheorem{proposition}[theorem]{Proposition}

\newcommand{\R}{\mathbb{R}}
\newcommand{\N}{\mathbb{N}}
\newcommand{\Params}{\mathcal{P}}
\newcommand{\grid}{\mathcal{T}_h}
\newcommand{\faces}{\mathcal{F}_h}
\newcommand{\Hdiv}{H_\textnormal{div}}

\begin{document}

\date{}

\title{A locally conservative reduced flux reconstruction for elliptic problems}
\author{Stephan Rave\footnotemark[1] \and Felix Schindler\footnotemark[1]}

\maketitle

\renewcommand{\thefootnote}{\fnsymbol{footnote}}
\footnotetext[1]{Mathematics Münster, Einsteinstr.~62, 48149 Münster, Germany, \url{{stephan.rave, felix.schindler}@uni-muenster.de}}

\begin{abstract}
  In the context of model order reduction of parametric elliptic problems, we present a methodology to reconstruct a conforming flux from a given reduced solution, that is locally conservative with respect to the underlying finite element grid.
  All components of the procedure depend separably on the parameter and allow for further use in offline/online decomposed computations, for instance in the context of a posterior error estimation or flow problems.
\end{abstract}

\section{Introduction}

In the context of model order reduction of elliptic problems \eqref{eq:weak_solution} by the reduced Basis (RB) method, it is sometimes desirable to obtain a locally conservative flux \eqref{eq:locally_conservative}, for instance for a posteriori error estimation \cite{OS2015} or in the context of flow problems \cite{KFH+2015}.
However, this is not straightforward for RB solutions \eqref{eq:RB_solution}, even if the underlying scheme \eqref{eq:DG_solution} yields locally conservative approximations.
While \cite{LR2018} ensure global mass conservation and \cite{CCS2018} ensure mass conservation with respect to few subdomains, we are not aware of any work ensuring local mass conservation with respect to the underlying grid in the RB context.
The present work closes this gap.
We will not give individual references from here on but refer to \cite{OS2015} and the references therein.

\section{Nonconforming approximation and diffusive flux reconstruction}

Given a bounded physical domain $\Omega \subset \R^d$ for $d =  2, 3$ with polygonal boundary and a compact parameter domain $\Params \subset \R^p$ for $p \in \N$, consider for each parameter $\mu \in \Params$ the weak solution $u_\mu \in H^1_0(\Omega) \subset H^1(\Omega) \subset L^2(\Omega)$ of
\begin{align}
  -\nabla\cdot\big( \sigma_\mu \nabla u_\mu \big) = f, &&\text{with } f \in L^2(\Omega) \text{ and parametric uniformly positive } \sigma_\mu \in L^\infty(\Omega)\text{ for all } \mu \in \Params,
  \label{eq:weak_solution}
\end{align}
the latter ensuring existence of a unique solution for each $\mu \in \Params$, owing to the Lax-Milgram Theorem.

Given a nonoverlapping conforming simplicial partition of $\Omega$ by a grid $\grid$ with disjoint elements $K \in \grid$, such that $\cup_{K \in \grid} \overline{K} = \overline{\Omega}$ and corresponding faces $\faces$, each $\Gamma \in \faces$ given by the intersection of exactly two distinct elements $K^-, K^+ \in \grid$, $\Gamma = \overline{K^-} \cap \overline{K^+}$ or by intersection of a single element $K^- \in \grid$ with the domain boundary, $\Gamma = \overline{K^-} \cap \partial \Omega$, we consider interior penalty (IP) discontinuous Galerkin (DG) approximations $u_{h, \mu} \in V_h := \{v \in L^2(\Omega) \,|\, v|_K \in \mathbb{P}_1(K)\;\forall K \in \grid\} \not\subset H^1_0(\Omega)$ of \eqref{eq:weak_solution}, where $\mathbb{P}_1(K)$ denotes the space of polynomials of degree at most 1 over a simplex $K \in \grid$, for $\mu \in \Params$ as the solution of
\begin{align}
  a_{h, \mu}(u_{h, \mu}, v_h) &\,= (f, v_h)_{L^2(\Omega)} \quad\text{for all } v_h \in V_h, \text{ with the symmetric IPDG bilinear form given by}
  \label{eq:DG_solution}\\
  a_{h, \mu}(u, v) &:= \sum_{K \in \grid} \int_K (\sigma_\mu \nabla u) \cdot \nabla v \,\textnormal{d}x
  + \sum_{\Gamma \in \faces} \int_\Gamma -\big<(\sigma_\mu \nabla v)\cdot n\big>[u] - \big<(\sigma_\mu \nabla u)\cdot n\big>[v] + \tfrac{\nu}{h_\Gamma}[u][v] \,\textnormal{d}s,
  \notag
\end{align}
with $h_\Gamma := \textnormal{diam}(\Gamma)$, the normal $n \in \R^d$ on each face pointing away from $K^-$ and, for a two-valued function $v \in V_h$, with the average and jump over a boundary face given by $\left<v\right> := [v] := v$ and on an inner face given by $\left<v\right> := \tfrac{1}{2}(v^- + v^+)$ and $[v] := v^- - v^+$, respectively, with $v^\pm := v|_{K^\pm}$.
As usual for IPDG methods, $a_{h, \mu}$ is continuous and coercive with respect to a suitable DG-norm on $V_h$ if the user-dependent penalty parameter $\nu > 0$ is chosen large enough.
Thus a unique IPDG solution exists for each $\mu \in \Params$, owing to the Lax-Milgram Theorem.

To obtain a conforming flux we invoke the diffusive flux reconstruction operator $R_h: V^h \to Q^h$, where $Q_h \subset \Hdiv(\Omega)$ denotes the $0$-st order Raviart-Thomas-N\'{e}d\'{e}lec space of vector-valued fluxes with continuous normal components with respect to $\faces$ and $\Hdiv(\Omega) := \{ q \in L^2(\Omega)^d \,|\, \nabla\cdot q \in L^2(\Omega) \}$.
We specify $R_{h, \mu}$ locally for some $v_h \in V_h$ by
\begin{align}
  \int_\Gamma R_{h, \mu}(v_h)\cdot n\,\textnormal{d}s = \int_\Gamma -\big<(\sigma_mu \nabla v_h)\cdot n\big> + \tfrac{\nu}{h_\Gamma}[v_h]\,\textnormal{d}s &&\text{for all } \Gamma \in \faces.
  \label{eq:flux_reconstruction}
\end{align}
If applied to the IPDG solution $u_{h \mu} \in V_h$, we obtain a conforming flux $t_{h, \mu} := R_{h, \mu}(u_{h, \mu}) \in \Hdiv(\grid, f)$ that is in addition \emph{locally conservative}, which we specify by the following set, with respect to $\grid$ and $g \in L^2(\Omega)$:
\begin{align}
  \Hdiv(\grid, g) := \big\{ q \in \Hdiv(\Omega) \;\big|\; (\nabla\cdot q, 1)_{L^2(K)} = (g, 1)_{L^2(K)}\quad \forall K \in \grid \big\}
  \label{eq:locally_conservative}
\end{align}
This local conservation property holds since the constant function which evaluates to 1 on $K \in \grid$ and to 0 elsewhere lies in $V_h$ for all $K \in \grid$, which can be used to show $R_{h, \mu}\, u_{h, \mu} \in \Hdiv(\grid, f)$ by applying Greens theorem, and by using \eqref{eq:flux_reconstruction} and \eqref{eq:DG_solution}.

\section{Model order reduction by reduced basis methods}

To obtain quickly available approximations of $u_{h, \mu}$ and $R_{h, \mu}(u_{h, \mu})$, for instance in the context of real-time or multi-query contexts, we seek a low-dimensional reduced space $V_n \subset V_h$ to approximate the set $\{u_{h, \mu} \,|\,\mu \in \Params\}$, with $n := \dim V_n \ll \dim V_h$, usually $n \approx 10^2$ while $\dim V_h \approx 10^{6 \text{ to } 9}$.
Presuming we are given $V_n$ (spanned by solution snapshots $u_{h, \mu}$ for certain $\mu \in \Params$ selected by a greedy algorithm, see Algorithm \ref{alg:greedy}), we obtain a reduced solution $u_{n, \mu} \in V_n$ for $\mu \in \Params$ as the unique solution of
\begin{align}
  a_{h, \mu}(u_{n, \mu}, v_n) = (f, v_n)_{L^2(\Omega)} \quad\text{for all } v_n \in V_n, \text{ namely the Galerkin projection of \eqref{eq:DG_solution} onto $V_n$.}
  \label{eq:RB_solution}
\end{align}
If $a_{h, \mu}$ depends separably on the parameter, that is $a_{h, \mu} = \sum_{\xi = 1}^\Xi \theta_\xi(\mu) a_{h, \xi}$ for parameter functionals $\theta_\xi: \Params \to \R$ and component bilinear forms $a_{h, \xi}$, we may precompute the projection of the $a_{h, \xi}$ onto $V_n$ once in an offline step and may later on, given $\mu \in \Params$, assemble and solve the linear system associated with \eqref{eq:RB_solution} with a complexity of $\mathcal{O}(\Xi n^2 + n^3)$, independent of $\dim V_h$.
To obtain quasi optimal $V_n$ we invoke a greedy algorithm which relies on an reliable a posteriori error estimate $\|u_{h, \mu} - u_{n, \mu}\| \leq \eta_n(\mu) := \underline{\sigma}_\mu \|\mathcal{R}_{h, \mu}(u_{n, \mu})\|_{-1}$, where $\underline{\sigma}_\mu > 0$ denotes the minimum of $\sigma_\mu$ over $\Omega$, and $\mathcal{R}_{h, \mu}$ denotes the residual associated with \eqref{eq:DG_solution}, the dual norm of which can be efficiently computed by means of its Riesz-representatives, using the separable decomposition of $a_{h, \mu}$.

\section{Locally conservative reduced flux reconstruction}

Unfortunately, the flux reconstruction of a reduced solution $u_{n, \mu}$ is in general not locally conservative, $R_{h, \mu}\, u_{n, \mu} \not\in \Hdiv(\grid, f)$, since the indicator functions for each $K \in \grid$ are missing in $V_n$.
To remedy the situation, we note that while $\Hdiv(\grid, f)$ is not a linear subspace of $\Hdiv(\Omega)$, $\Hdiv(\grid, 0) \subset \Hdiv(\Omega)$ is and since $\Hdiv(\grid, f) = \Hdiv(\grid, 0) + t_f$ for $t_f \in \Hdiv(\grid, f)$ we build a reduced flux space $Q_n \subset \Hdiv(\grid, 0)$ along $V_n$ (lines 1, 2, 3 in Algorithm \ref{alg:greedy}).
\begin{algorithm2e}[h]
\DontPrintSemicolon
\KwIn{finite dimensional training set $\Params_\textnormal{train} \subset \Params$, tolerance $\delta > 0$}
\KwOut{reduced spaces $V_0 \subset \dots \subset V_n$, reduced flux spaces $Q_0 \subset \dots \subset Q_n$, $t_f \in \Hdiv(\grid, f)$}
\SetKwFunction{ErrEst}{ErrEst}
\SetKwFunction{RBSolve}{RBSolve}
\SetKwFunction{Solve}{Solve}
\SetKwFunction{Reduce}{Reduce}
$n \leftarrow 0$,\ \ $V_0 \leftarrow \{0\}$,\ \ $Q_0 \leftarrow \{0\}$\;
\nl$t_f := R_{h, 1}\,u_{h, f} \in \Hdiv(\grid, f)$ with $u_{h, f} \in V_h$ solution of \eqref{eq:DG_solution} with $\sigma_\mu = 1$\;
\While{$\max_{\mu \in \Params_\textnormal{train}} \eta_n(\mu) > \delta$}{
  $n \leftarrow n + 1$,\ \ $\mu^* \leftarrow \operatorname{argmax}_{\mu \in \Params_\textnormal{train}} \eta_n(\mu)$\;
  $u_{h, \mu^*} \in V_h \leftarrow $ solution of \eqref{eq:DG_solution} for $\mu^*$\;
  $V_n \leftarrow \operatorname{span}(V_{n-1} \cup \{u_{\mu^*}\} )$\;
  \nl$t_{h, \mu^*} := R_{h, \mu^*}\, u_{h, \mu^*} \in \Hdiv(\grid, f)$ locally conservative flux reconstruction \eqref{eq:flux_reconstruction} for $\mu^*$\;
  \nl$Q_n \leftarrow \operatorname{span}(Q_{n-1} \cup \{t_{h, \mu^*, 0}\})$ with $t_{h, \mu^*, 0} := t_{h, \mu} - t_f \in \Hdiv(\grid, 0)$\;
}
\caption{weak discrete greedy basis generation with simultaneous locally conservative flux reconstruction}
\label{alg:greedy}
\end{algorithm2e}
We then define a reduced projection $\Pi_0: V_n \to Q_n$ for $u_{n, \mu} \in V_n$ as the solution of $(\Pi_0\, v_n, q_n)_{\Hdiv(\Omega)} = (R_{h, \mu}\, u_{n, \mu} - t_f, q_n)_{\Hdiv(\Omega)}$ for all $q_n \in Q_n$.
We thus obtain the following result, which is evident by construction.
\begin{proposition}
  Given a reduced solution $u_{n, \mu} \in V_n$ for $\mu \in \Params$, it holds that the reduced flux reconstruction $t_{n, \mu} := \Pi_0\,u_{n, \mu} + t_f$ is locally conservative, that is $t_{n, \mu} \in \Hdiv(\grid, f)$, and all required computations depend separably on $\mu$.
\end{proposition}

Extensions to higher orders, other grids and other IPDG schemes are straighforward, following \cite{OS2015}.
Error control of the flux reconstruction is subject to future work.

\section*{Acknowledgement}

  This work was funded by the Deutsche Forschungsgemeinschaft (DFG, German Research Foundation) under Germany ’s Excellence Strategy – EXC 2044 – 390685587, Mathematics Münster: Dynamics – Geometry - Structure,
  by the DFG under contract SCHI 1493/1-1 and by the German Federal Ministry of Education and Research (BMBF) under contract 05M18PMA.

\providecommand{\WileyBibTextsc}{}
\let\textsc\WileyBibTextsc
\providecommand{\othercit}{}
\providecommand{\jr}[1]{#1}
\providecommand{\etal}{~et~al.}

\end{document}